\documentclass[12pt]{article}
\usepackage{amsmath, amsfonts,amsthm}
\usepackage{latexsym}
\usepackage{graphicx}
\usepackage{color}
\usepackage{comment}

\usepackage[a4paper, total={7in, 9in}]{geometry}

\title{Electrostatic Origins of the Dirichlet Principle}

\author{Steven Deckelman\footnote{Department of Mathematics, Statistics and Computer Science,
University of Wisconsin-Stout,\textsf{ deckelmans@uwstout.edu}} }

\date{\today}

\newtheorem{mydef}{Definition}
\newtheorem*{cthm*}{Contingent Theorem}

\newtheorem*{thm*}{Theorem}
\newtheorem*{def*}{Definition}
\newtheorem{post}{Postulate}
\newtheorem{remark}{Remark}

\providecommand{\keywords}[1]
{
  \small	
  \textbf{\textit{Keywords---}} #1
}

\begin{document}
\maketitle

\begin{abstract}
The Dirichlet Principle is an approach to solving the Dirichlet problem by means of a
Dirichlet energy integral. It is part of the folklore of mathematics that the genesis of this argument was motivated by physical analogy involving electrostatic fields. The story goes something like this: 
If an electrostatic potential is prescribed on the boundary of a region, it will extend to a potential in the interior of the region which is harmonic when the electric field is in stable equilibrium, and that electrostatic field has minimum Dirichlet energy. The details of this argument are seldom given and where they are, they are typically scant, redacted, and speculative while often omitting either physics details or mathematics details. The purpose of  this article is to  give a detailed reconstruction of the electrostatic argument by combining accounts in several contemporary and historical disparate sources. Particular attention is given to explaining the frequently omitted physics and mathematical details and how they fit together to give the physical motivation.
\end{abstract}

\


\newpage

\section{Prologue}

\begin{quotation}
{\it Dirichlet published several papers on physical questions as 
well as on problems of pure mathematics. An interest and a creative 
power in both subjects seems nowadays to be rare, but it was not uncommon in those years. }- A.F. Monna
\end{quotation}

Peter Gustav Lejeune Dirichlet, 1805-1859, was one of the great luminaries of the German renaissance of mathematics in the nineteenth century. The successor of Gauss at G\"{o}ttingen and predecessor of  Riemann, Dirichlet remains a fascinating exemplar of what Constance Reid called G\"{o}ttingen's  mathematical-physical tradition in \cite{Reid}. In 1825 Dirichlet was catapulted to fame as a twenty-year-old student without a degree when he proved Fermat's Last Theorem for $n=5$. Unable to get a traditional doctorate because his lack of fluency in spoken Latin\footnote{Dirichlet did have a considerable command of written Latin as is evident from his collected works.} prevented him from participating in the required disputation of his thesis, he was awarded an honorary doctorate instead.\footnote{He did later earn his habilitation.} A cigar smoker, Dirichlet suffered a devastating heart attack in the summer 1858 and died in 1859. See \cite{Els} and \cite{MB}.

This article deals with certain questions from the history of mathematics that touch on the relationship between mathematics and physics in the nineteenth century.
The {\it Dirichlet principle} is a method for solving the Dirichlet problem by minimizing the Dirichlet integral. If $\Omega\subset\mathbb{R}^3$ this integral is given below in equation
\ref{ei}. The minimization problem is over all sufficiently smooth extensions $U$ of a prescribed boundary function $f$. This same integral is called an energy integral and purportedly is heuristically justified by physical analogy with electrostatic fields. The provenance of this physical argument is obscure and in the relatively few places in the literature
which do mention it, details are typically scant and highly redacted with either omissions of mathematical  or physical details or both. In what follows it will be argued that the reason for this state of affairs is that electrostatic argument was part of the mathematical folklore of the nineteenth century and not written down in any detailed form. We will also argue that notwithstanding this, there was an actual electrostatic argument that can be reconstructed, and we give a detailed reconstruction. We do this by elaborating upon and filling in the details of  three redacted accounts given by Franz Grube, Oliver Dimon Kellogg, and A.F. Monna. The  evidence for underlying historical thesis of our account is based on a number of observations by prominent mathematicians that testify to the uncertain and indeed unknown form of the original physical argument. These are expressed through various quotations that appear  below.

The crux of this article is contained in 
sections \ref{ho}, \ref{dcalc}, \ref{eim}, and \ref{cv} below.  Although this article is intended to be an expository/survey, some of the details appearing in section \ref{dcalc} are new in the sense that we are unaware of any place in the literature where they can be found. This is the basis for using the term reconstruction for the argument to calculate the Dirichlet energy. In \cite{Kel} Kellogg speculates about the general form the original argument may have taken on page 278 of chapter XI, gives the equation \ref{energysummary}  below and the energy integral that results from it in equation \ref{complete-energy} below, but Kellogg omits the proof or details of the calculation. Perhaps he considered it obvious. 
 We provide these below. Indeed, we are not completely certain the arguments we give are what Kellogg had in mind; we offer our derivation as a plausible possibility.
A ubiquitous issue in any historical reconstruction is the trade-off between historical fidelity of capturing the way people at the time conceptualized a subject versus  intelligibility for contemporary readers enhanced by using modern language and constructs. The electron, for example, was not even discovered until 1897, some thirty-eight years after Dirichlet's death. Since we hope this article will be of interest to a general audience, we take the latter approach and freely make use of modern constructs when we believe they contribute to intelligibility.

\section{Physical Mathematics, the Dirichlet Principle, and Analogical Reasoning}

\begin{quotation}
{\it The science of physics does not only give us (mathematicians) an opportunity to solve problems, but helps us to discover the means of solving them, and it does this in two ways: it leads us to anticipate the solution and suggests suitable lines of argument.} -  Henri Poincar\'{e}
\end{quotation}

The word {\it physical mathematics}, in contradistinction to mathematical physics, has been used to describe a thought process in mathematics guided  by physical interpretations or analogies.  Felix Klein called this the physical method {\it (Die physikalische Methode)} in \cite{K}, and argued that it strongly influenced Riemann's thinking about algebraic integrals and Riemann surfaces. In particular Klein cited Riemann's use of the Dirichlet principle as evidence of physical thinking in Riemann's thought. There are many historical examples of physical mathematics including Archimedes' mechanical approach to the quadrature of the parabola \cite{A1},  Johann Bernoulli's ``solution'' of the Brachistochrone problem  using optics and Snell's law \cite{P}, the Gauss-Lucas theorem via electrostatics (or fluid dynamics) \cite{Mar} and even an approach to the Riemann hypothesis via quantum mechanics (quantum Hamiltonians), the Hilbert-Polya conjecture \cite{B}. In the following pages we survey  and reconstruct in detail one of the most celebrated examples in the history of mathematics, the classical Dirichlet principle and the physical reasoning often ascribed to it.  Only elementary calculus and simple ideas about electricity and electrostatic fields are used.

By the classical Dirichlet problem we mean the boundary value problem for the Laplacian.
 Let $\Omega$ be a bounded, nonempty, connected, open set in $\mathbb{R}^n$ and let be given a ``sufficiently smooth'' real valued function $f$ on the boundary $S:=\partial\Omega$, where $S$ is also assumed to be ``sufficiently smooth''. In Dirichlet's day the open set 
$\Omega\subset \mathbb{R}^3$ could be interpreted physically as a conductor, with boundary corresponding to the surface of the conductor.  A problem of interest at the time was that of determining what today we'd call the electric field determined by a prescribed potential or charge distribution in a conductor. In mathematical terms the problem is to 
find a  function $u$ (or at least show that one exists), defined on the closure $\overline{\Omega}$ of $\Omega$, harmonic in $\Omega$ and such that $u=f$ on $S$ for some prescribed function $f$ . By harmonic in $\mathbb{R}^3$, we mean a twice continuously differentiable function $u(x,y,z)$ satisfying $\Delta u=0$ in $\Omega$, with $u$ continuous up to the boundary.  Here 
$$
\Delta=\frac{\partial^2}{\partial x^2}+\frac{\partial^2}{\partial y^2}+\frac{\partial^2}{\partial z^2}.
$$
That is,
\begin{equation}\label{bvp1}
\Delta u=0 \text{ on }\Omega,\qquad u=f \text{ on } S.
\end{equation}

The {\it Dirichlet Principle} is the moniker originated  by Bernhard Riemann in his 1857 paper on Abelian functions \cite{R2}, \cite{RT}, to describe P.G.L. Dirichlet's method for solving the boundary value problem for the Laplacian (the Dirichlet problem) by means of an
 energy integral. He had made use of it earlier in his 1851 inaugural dissertation \cite{R1} without reference to Dirichlet. Although the history of mathematics is replete with many instances of misattribution, Riemann did write of a  ``principle of Dirichlet''.
Possibly an instance of the Boyer-Stigler law of mathematical eponymy, it has been suggested by Felix Klein that Riemann coined this moniker owing more to the veneration he held for his teacher Dirichlet \cite{Pap}, than as being the source of the original idea, which can be found also in the writings of Green, Gauss, and Willam Thomson.  The word energy is suggestive of a physical interpretation and 
part of the folklore of mathematics is that Dirichlet's original chain of thought was guided by physical considerations involving the energy of electrostatic fields. 
We use the term folklore because, possibly as a consequence of it not being mathematically rigorous, it doesn't appear to be 
formally recorded in the literature, and where it is, it is typically in the form of a speculative sketch of the argument. See for example \cite{Gr}, \cite{Kel}, \cite{M}.
 Nonetheless we argue that this historical example of reasoning by physical analogy remains of interest because it illustrates the way physics  influenced the creative work of mathematicians of the past.

Physical analogy is a particular type of reasoning by analogy. Reasoning by analogy is by definition non-rigorous. Its conclusions cannot be substituted for mathematical proof. But it can guide reasoning in directions that can be placed on more firm mathematical foundations with the right tools. The flip side of this is that analogical reasoning does give us license to take steps that we may know (or believe) are logically unjustified but consonant with existing analogies. We employ this license freely below.

\section{The Dirichlet Principle}

\begin{quotation}
{\it
 In 
many cases, notably integrals of algebraic functions and their inverse functions, there is a principle which Dirichlet employed for this purpose. Probably 
inspired by a similar idea of Gauss, he applied this principle to a function of 
three variables satisfying Laplace's partial differential equation in his lectures 
over several years on forces obeying an inverse square law \dots
 Perhaps I may refer, in regard to some 
related researches, to the account in my doctoral dissertation. \footnote{Riemann's doctoral thesis, Foundations 
for a general theory of functions of a complex variable, G\"{o}ttingen 1851.}
} -G. B. Riemann
\end{quotation}

 In this section we describe in greater detail Dirichlet's approach to the boundary value problem \ref{bvp1}.  We will describe the problem in space, since this is the setting that  Dirichlet worked in. The reader may find it notable the physical heuristic derivation presented below is tightly tethered to three-dimensional space and is not adapted to the plane.\footnote{A heuristic argument more directly adapted to the plane could probably be formulated by replacing Newtonian potentials in space by logarithmic potentials in the plane.} Once the correct functional to be minimized, however, is identified and given in equation \ref{ei} below, the more rigorous variational argument detailed in section \ref{cv} of this article works in any dimension.

{\bf Dirichlet Principle:} Let the prescribed boundary function $f$ be given. 
By the class of  admissible functions, which needs yet to be defined, but for which we will provisionally assume
is meant any twice continuously differentiable function $U$ defined in $\Omega$ and continuous on 
$\overline{\Omega}$\footnote{The topological closure of $\Omega$.} with the same being true of  its exterior $\mathbb{R}^3\setminus \overline{\Omega}$   that is equal to $f$ on $S$. Then there will exist an admissible function $u$ that minimizes the integral\footnote{In our derivation of the Dirichlet energy integral below we will express the energy as an integral over all of $\mathbb{R}^3$ but ultimately confine attention to integration over $\Omega$ for the interior Dirichlet problem.}
\begin{equation} \label{ei}
D(U)=\iiint_\Omega \left[\left( \frac{\partial U}{\partial x}\right)^2 +\left( \frac{\partial U}{\partial y}\right)^2+\left( \frac{\partial U}{\partial z}\right)^2\right] dV
\end{equation}
and this $u$ satisfies the boundary value problem \ref{bvp1}. Here $dV$ is the volume measure $dxdydz$. We will explain the sense in which this iconic integral, so often quoted without derivation, represents energy in equation \ref{complete-energy} below.

Note that the admissible functions $U$ appearing in the integral \ref{ei} need not be harmonic.  The integral $D(U)$ is often called the energy integral because it can be shown to represent the total potential energy of a mass distribution (charge distribution) in space associated with the vector field for which $U$ is potential, $F=\nabla U$ (or defined by $-F=\nabla U$ in physics texts)
 off $S$ (both inside and outside of $\Omega$) arising from an inverse square law of attraction or repulsion. 
In this article we use the symbol $F$ to represent an electric field rather than the more customary $E$, because we reserve the latter symbol for energy.

 Conceptually, the physical argument might be described as follows. By the principle of potential energy minimization, the mass distribution need not be in stable equilibrium, but we imagine that the mass distribution ``seeks out'' a stable equillibrium. We envision a given mass distribution in space as able to rearrange itself, subject to a prescribed boundary potential that would be implemented through an appropriate experimental set up,  into a configuration of minimum potential energy. That the mass distribution associated with an arbitary potential $U$ can be recovered from $U$ is based on an old theorem of Poisson. In what follows these ideas will become clearer.

\section{Law of Conservation of Energy; Work as Energy} \label{sec:WE}

\begin{quotation}
{\it The origin of the principle is in physics, although this is no longer apparent from the modern formulation. Originally the integral that was to be minimized meant something like the energy of a system $\dots$  Hence it is not surprising that these mathematicians frequently used physical arguments, though only with heuristic value.}
-A.F. Monna
\end{quotation}

Suppose we have a vector field in a region of space, say $\mathbb{R}^3$. In physics energy is defined as the ability to do work. If a particle moves along a path $x=x(t)$,
$y=y(t)$, $z=z(t)$, $a\leq t\leq b$, acted upon by a force $F(\mathbf{X})$,
then the work $W$  done will be given by the line integral
$$
W=\int_C F(\mathbf{X})\cdot d\mathbf{X},
$$
where $C$ is the path in space and $\mathbf{X}=\mathbf{X}(t)=(x(t),y(t),z(t))$. 

From calculus we know that these line integrals will be independent of path precisely when the vector field $F$ is conservative. In such a case we would write
$$
W=\int_P^Q F(\mathbf{X})\cdot d\mathbf{X},
$$
where $P=\mathbf{X}(a)$ and $Q=\mathbf{X}(b)$.

Furthermore, the vector field will be conservative precisely when it has an antigradient (antiderivative)-a scalar function $U(x,y,z)$ for which
\begin{equation}\label{ag}
 \nabla U=F .
\end{equation}
Moreover in this case there will be a version of the fundamental theorem of calculus, namely,
\begin{equation}  \label{work1}
U(Q)-U(P)=\int_P^Q \nabla U(\mathbf{X})\cdot d\mathbf{X}=\int_P^Q F(\mathbf{X})\cdot d\mathbf{X}=W.
\end{equation}
The antigradient $U$ is also called the {\it potential function} for the vector field, or simply the {\it potential}. Note that equation \ref{work1} says that the work done by a particle moving from point $P$ to point $Q$ in a conservative field is the difference of the antigradient function $U(Q)-U(P)$. If $Q$ is a fixed reference point, $U(x,y,z)$ is, by definition, the potential energy. Note also that the {\it potential energy} (being the negative of an antigradient) is defined only up to an arbitrary constant.  The potential energy is really a property of the vector field. This is why it is sometimes called the energy of position. It probably would be more accurate to call it the energy of the configuration of the forces in the vector field. Note as well that potential energy is only defined for conservative vector fields.

Now lets take a slightly different tack. If we apply Newton's Second Law to $F$,
$F=ma$, with $m$ the mass of our particle\footnote{In what follows we shall use the terms ``mass'', ``charge'', ``distributions'', and ``densities'' interchangeably. We think of charge as ``electrical mass''.} , we can write our force as
$$
F=m\left(\frac{d^2x}{dt^2},\frac{d^2y}{dt^2},\frac{d^2z}{dt^2}\right).
$$
We can then write the work done between $P$ and $Q$ as the line integral
$$
W=\int_P^Q m\left(\frac{d^2x}{dt^2},\frac{d^2y}{dt^2},\frac{d^2z}{dt^2}\right)\cdot d\mathbf{X}
=m\int_P^Q \left(\frac{d^2x}{dt^2}dx+\frac{d^2y}{dt^2}dy+\frac{d^2z}{dt^2}dz\right)
$$

$$
=m\int_P^Q \left(\frac{d^2x}{dt^2}\frac{dx}{dt}+\frac{d^2y}{dt^2}\frac{dy}{dt}+\frac{d^2z}{dt^2}\frac{dz}{dt}\right)dt
=\frac{1}{2}m\int_P^Q\frac{d}{dt}\left[ \left(\frac{dx}{dt}\right)^2 +\left(\frac{dy}{dt}\right)^2+\left(\frac{dz}{dt}\right)^2\right]dt
$$

$$
=\frac{1}{2}m\int_P^Q \frac{d}{dt}\left\vert \mathbf{V}(t)\right\vert^2 dt=\frac{1}{2}m\left(\left\vert \mathbf{V}(Q)\right\vert^2-\left\vert \mathbf{V}(P)\right\vert^2\right).
$$
Here $\mathbf{V}$ is the velocity vector $\mathbf{V}=\left(\frac{dx}{dt},\frac{dy}{dt},\frac{dz}{dt}\right)$. The quantity $\frac{1}{2}m\left\vert \mathbf{V}\right\vert^2$
is by definition called the {\it kinetic energy}. Comparing this with equation \ref{work1} we see we have two expressions for the work done between points $P$ and $Q$:
\begin{equation} \label{preminus}
W=U(Q)-U(P)=\frac{1}{2}m\left(\left\vert \mathbf{V}(Q)\right\vert^2-\left\vert \mathbf{V}(P)\right\vert^2\right).
\end{equation}
This formula assumes a nicer form if we redefine the function $U$ (following the physics convention) in equation \ref{ag} as
$$
\nabla U=-F.
$$
Equation \ref{preminus} then takes the form
\begin{equation} \label{LCE}
W_{PQ}=U(P)-U(Q)=\frac{1}{2}m\left(\left\vert \mathbf{V}(Q)\right\vert^2-\left\vert \mathbf{V}(P)\right\vert^2\right).
\end{equation}

In words we may say work done is equal to the  change in kinetic energy or to the negative of the
difference or change in potential energy. The law \ref{LCE} is called the {\it Law of Conservation of Energy} and usually written as
\begin{equation}\label{lce}
U(P)+\frac{1}{2}m\left\vert \mathbf{V}(P)\right\vert^2=U(Q)+\frac{1}{2}m\left\vert \mathbf{V}(Q)\right\vert^2.
\end{equation}

When our field is conservative, the law of conservation of energy \ref{LCE} says that the work done is the (negative of) the change in potential energy, so that
\begin{equation}\label{workenergy}
U(P)=U(Q)+W_{PQ}.
\end{equation}
Since potentials are defined only up to an additive constant, equation \ref{workenergy} justifies identifying work done with potential energy.\footnote{Conceptually work and energy are different but closely related. One Joule of work is done when a one Newton force moves an object one meter. One Joule of energy is the energy capable of applying a one Newton force in moving said object one meter by transmission of energy to the object through motion.}
 Incidentally, the quantity appearing in equation \ref{lce},

\begin{equation} \label{tme}
E=U+\frac{1}{2}m\left\vert \mathbf{V}\right\vert^2
\end{equation}
is called the total mechanical energy.

\section{Potentials in Mathematics and Physics}\label{3p}

\begin{quotation}
{\it
The word "potential" is a hideous misnomer because it inevitably reminds you of potential energy. This is particularly confusing, because there is a connection between "potential" and "potential energy," as you will see in Section
2.4. I'm sorry that it is impossible to escape this word. The best I can do is to insist once and for all that "potential" and "potential energy" are completely different
terms and should, by all rights, have completely different names. (Incidentally, a surface over which the potential is constant is called an equipotential.)
- David J. Griffiths}
\end{quotation}

In physics the term {\it potential} can refer to several interrelated but different concepts.
There are different kinds of potentials. This can be  a source of confusion as  reflected in the quotation so dolefully lamented in the  
   classic undergraduate text on electrodynamics \cite{G}.
Perhaps Professor Griffith's position is a little extreme here, but there are some nuanced differences important to be aware of. At a very general level, mathematicians often think of a potential as a convolution of function and a measure. 
This general notion can take myriad forms.   In the context we are considering the  term ``potential'' can refer to any of the following:

\begin{itemize}
\item The potential function of a vector field. This is a scalar function  $U$ which acts as an antigradient or antiderivative of a field $F$. In other words $\nabla U=F$, or what is usually more convenient,
$\nabla U=-F$. The existence of a potential function means the field is conservative. The potential function $U$ need NOT be harmonic.
\item The potential energy of the field $F$. This will be given by the potential function $U$ and appears in the field's law of conservation of energy \ref{lce} and the expression for total mechanical energy \ref{tme}.
\item The volume (or spatial) potential of a charge distribution.  This is a function.
\item The surface potential of a charge distribution. This is a function.
\item The total potential of a charge distribution, the sum of the volume and surface potentials.\footnote{The spatial, surface, and total potentials will be explained below.}
\item  The potential energy of a mass/charge distribution in itself. This is a scalar. That is, the work/energy needed to assemble a distribution of mass from some standard distribution or location such as infinity. In some older books like Kellogg \cite{Kel} this is sometimes called a {\it self-potential}.
\item The self-potential of a volume distribution. This is a scalar. It represents the energy needed to assemble the spatial distribution.
\item The self-potential of a surface distribution. This is a scalar. It represents the energy needed to assemble the surface distribution.
\item The mutual potential of a system consisting of both volume and surface parts. This is a scalar.  It is an interaction energy that arises between the volume and surface parts of a distribution. The self-potential of the total system is the sum of the volume and surface self-potentials as well as the mutual potential.
\item A harmonic function. Vector fields of charge distributions in stable equilibrium will have harmonic potential functions.
\end{itemize}

In electrostatics, physicists distinguish between the concepts of (electrostatic) potential energy, the electric field, and the electric potential. Understanding these distinctions may be helpful in understanding the physical argument behind the Dirichlet principle.
As in section \ref{sec:WE}, $F$ and $U$ denote a vector field and it's potential, i.e.
$$
\nabla U=-F.
$$
In electrostatics the mutual force of attraction or repulsion between charges at rest is governed by Coulomb's law. Suppose $Q$ represents a fixed electric charge at some point in space. The fixed charge $Q$ creates a ``potential''  electric field, that comes into existence through interaction with other charges.
If $q$ is a ``test charge'', 
at some location, Coulomb's law will describe the electric field created by  $Q$ acting on $q$. This is what is designated by  equations \ref{coulomblaw}, and \ref{coulomb} below.
The electrostatic force $F$ appearing in Coulomb's formula is measured in Newtons. 
Work is measured in Joules (equivalently Newton-Meters).  The potential energy $U$ is also measured in Joules. Strictly speaking, physicists would not call $F$ the electric field, but would make a more careful distinction between the ``acting charge'' $Q$ and the ``acted upon charge'' $q$. This more careful definition of an electric field is then defined by\footnote{Note that while we here use $\mathcal{E}$ for the electric field, this symbol is often used for electromotive force (emf), a scalar, in many current physics texts. }
$$
\mathcal{E}=\frac{F}{q}.
$$
The units of $\mathcal{E}$ will then be Newtons per Coulomb. Dividing both sides of equation \ref{work1} by $q$ yields
$$
\frac{U(Q)-U(P)}{q}=\int_P^Q \frac{\nabla U(\mathbf{X})}{q}\cdot d\mathbf{X}=\int_P^Q \frac{ F(\mathbf{X})}{q}\cdot d\mathbf{X}.
$$
Physicists then define a new scaler field
$$
V=\frac{U}{q},
$$
whose units are thus Joules per Coulomb (equivalently Volts). $V$ is then called the electric potential since it will be the antigradient of $\mathcal{E}$, that is
$$
\nabla V=-\mathcal{E}.
$$
The function $V$ is also called the electric potential, which is different from the electrostatic potential (energy) $U$.  Again, $V$ is measured in Volts while $U$ is measured in Joules. The units of $\mathcal{E}$ are measured in Newtons per Coulomb  (equivalently Volts per Meter).  

We observe that if $q=1$, then $V=U$ and $\mathcal{E}=F$. For the remainder of this article we will assume our test charge $q=1$ and refer to  $U$ and $F$ as the electric potential and electric field. But it should be kept in mind that although the mathematical quantities will be the same, the physical conceptual meanings as testified by the different units will be different.

In what follows we will describe how the electrostatic argument implies
 the relationship between these four sundry notions of  potential are as follows.
\begin{eqnarray} 
{\textrm Charge \  Distribution }\  \rho &\longleftrightarrow\text{ Potential Function }U(x,y,z) \nonumber\\
&\longleftrightarrow\text{ Vector Field }F. \label{bd}
\end{eqnarray}
Moreover, when the charge distribution $\rho$ is in stable equilibrium $\longleftrightarrow U(x,y,z)$ is harmonic.\footnote{We are assuming the stable equilibrium is unique.}

\section{The Principle of Potential Energy Minimization}

\begin{quotation}
{\it Riemann, as we know, used
Dirichlet's Principle in their place in his writings. But I have
no doubt that he started from precisely those physical problems,
and then, in order to give what was physically evident the
support of mathematical reasoning, he afterwards substituted
Dirichlet's Principle. }
-Felix Klein
\end{quotation}

The principle of potential energy minimization is the principle that a system in stable equilibrium at rest corresponds to minimum potential energy.
It seems to have been first stated in 1798 by J.L. Lagrange  in his famous
{\it M\'{e}canique Analytique} and ``proved'' \footnote{The principle of potential energy minimization may be as much a philosophical principle as it is mathematical principle. It appears to be an idealization of the observation that unstable rest states cannot practically persist in nature.}
 later by P.G. Lejeune Dirichlet in \cite{PGLD}.  Dirichlet himself pointed out that  the basic idea goes back to Daniel Bernoulli.
A statement is given in Brauer and Nohel's book \cite{BN}:
\begin{post}\label{postulate}
If in a certain rest position a conservative mechanical system has minimum potential energy, then this position corresponds to a stable equilibrium;
if the rest position does not correspond to  minimum potential energy, then the equilibrium is unstable.
\end{post}

By an electrical system we will mean a set of charge configurations defined by some set of physical constraints.
We will apply the minimum potential energy principle to these electrical systems and assume the following:
\begin{post} \label{es}
If in a certain rest position a conservative electrical system has minimum potential energy, then this system corresponds to a stable equilibrium; if the rest position does not correspond to minimum potential energy, then the equilibrium is unstable. 
\end{post}
For additional historical details on Dirichlet's ideas in this area, see Uta Merzbach's recent biography \cite{MB}.

Hamiltonian mechanics furnishes a particularly transparent means of proving this principle.
See for example \cite{CH}, pp 242-244, as well as  \cite{U}, \cite{BN}, and chapter one of \cite{Ma}.

\section{The Basic Datum: Potentials, Charge Distributions, and Electric Fields}

\begin{quotation}
{\it In the middle of the twentieth century it was attempted to divide physics and mathematics. The consequences
turned out to be catastrophic. Whole generations of mathematicians grew up without knowing half of their
science and, of course, in total ignorance of any other sciences. They first began teaching their ugly scholastic
pseudo-mathematics to their students, then to schoolchildren (forgetting Hardy's warning that ugly mathematics
has no permanent place under the Sun).}-V.I. Arnold
\end{quotation}

In order to make the classical  electrostatic argument as clear as possible to modern readers, we introduce some notation and definitions in this section.

We imagine a bounded conductor or cavity in space and model it by a domain (non-empty open connected set),
$\Omega\subset\mathbb{R}^3$ having boundary $S=\partial\Omega$. We then  decompose    
$\mathbb{R}^3 $	 into the union $\Omega\cup S\cup ext(\Omega)$. 
Here $ext(\Omega)$ is the topological exterior of $\Omega$.

\begin{mydef}\label{totaldensity}
 By a charge distribution we mean a mass distribution of negative 
charge\footnote{We choose negative rather than positive charge to emphasize a conceptual model in which electrons  or anions would be analogous to particles of electric charge. Of course, conventional (positive) charge could be used as well.}
 with both a spatial component on $\Omega$ and $ext(\Omega)$ as well as a surface component on $S$. Each such charge distribution can be identified with a spatial density function $\rho_1(x,y,z)$ and a surface density function $\rho_2(\zeta)$, $\zeta\in S$. By a slight abuse of notation we will use the symbol $\rho$ to represent either of these densities when it is clear by context. 
\end{mydef}
The idea here is that we want to think of distributions of charge both in space as well as on the surface of the conductor. By the measure/mass or total charge of a distribution $\rho$ we mean the sum of these. More specifically, if our charge distribution is supported on $E=F\cup G$ where $F\subset S^c$ and $G\subset S$, then the total charge would be calculated as
$$
\iiint_{F}\rho_1(x,y,z)dV +\iint_G \rho_2(\zeta)dS(\zeta).
$$

We emphasize that total charge or electrical mass (measure) of the charge distribution given in  definition \ref{totaldensity}  consists of the sums of the spatial and surface charges (measures).

\begin{mydef} \label{datum}
By the basic datum we mean the triple $(U,\rho,F)$ referenced above in equation \ref{bd} as the potential function $U$ for the  electric field $F$ of the charge distribution $\rho$, possibly satisfying some electrical system constraint such as a boundary condition.
 We will sometimes use these interchangeably or collectively. 
\end{mydef}

A basic datum $(U,\rho,F)$ is the mathematical representation of  an ``electrical system'' such as the one  described  in postulate \ref{es}.

\begin{mydef} \label{admiss}
Given a continuous function $f:S\to \mathbb{R}$, denote by $\mathcal{A}_f$ the class of all twice continuously differentiable extensions $U$ of $f$ to $S^c$ continuous up to the boundary $S$ in
 $\mathbb{R}^3$. We call $\mathcal{A}_f$ the class of admissible extensions of $f$.
\footnote{Functions $U:\mathbb{R}^3\to \mathbb{R}$ that continuously extend $f$ to a twice continuously differentiable function on $\Omega$ and a twice continuously differentiable function on $\Omega^c$. Notice $U$ need not be  twice continuously differentiable on all of $\mathbb{R}^3$. }
\end{mydef}

\section{Heuristic Origins of the Dirichlet Energy Integral}\label{ho}

\begin{quotation}
{\it  I am acutely aware of the fact that the marriage between mathematics and physics, which was so enormously fruitful in past centuries, has recently ended in divorce.}- Freeman Dyson
\end{quotation}

To show where the Dirichlet energy integral comes from, in this section we follow the ``classical'' approach 
of first deriving formulae for discrete mass distributions and then generalize this by analogy to  continuous mass distributions.  This was the technique used by Grube, Kellogg, and Monna, but our exposition most closely follows Kellogg's.

Recall that one of the notions of potential energy is the energy of a mass distribution itself.  Suppose we consider a law of attraction or repulsion like Coulomb's law. Such fields are conservative. The magnitude of a force between two charges at rest is directly proportional to the product of the charges and inversely proportional to the square of the distance separating them.
In physics books this is often written in vector form using SI units as
\begin{equation} \label{coulomblaw}
F=\frac{1}{4\pi\epsilon_0}\frac{q_1 q_2}{r^2}\vec{\mathbf{e}},
\end{equation}
where $\epsilon_0$ is a physical constant called the permittivity of space, $q_1$ and $q_2$ are the two charges in Coulombs and $\vec{\mathbf{e}}$ is a unit vector between the charge locations whose direction depends on the polarity (positive or negative) of the charges, i.e. whether attraction or repulsion is involved. For simplicity we assume units are chosen so that the magnitude of the force may be written as
\begin{equation}\label{coulomb}
|F|=\frac{m_1 m_2}{r_{12}^2}.
\end{equation}
We think of the $m_1$ and $m_2$ as ``electrical mass'' of the charges to emphasize the analogy between electric charge and matter.  The symbol $r_{ij}$ represents the distance  between two charges $m_i$ and $m_j$.

Suppose an assemblage of electrical masses $m_1$, $m_2$,...,$m_n$  is formed at positions
$P_1$, $P_2$,...,$P_n$ in space by  moving each $m_i$ separately, and sequentially, to position $P_i$ from infinity.  What will be the work done  forming this assemblage or configuration of mass? 
For conceptual simplicity, assume for the moment our mass system exists in  a plane. We recall the classical derivation. See for example \cite{Feyn}. For simplicity let us assume all charges are negative\footnote{We may think of these as electrons. This is somewhat at variance with conventional charge which is usually thought of as positively charged particles in textbooks.} so that the force between them is one of repulsion. The work done in moving $m_1$ to $P_1$ is zero because there are initially no charges in space to act on $m_1$. Once $m_1$ occupies $P_1$, Coulomb's law establishes an electrostatic field  according to equation \ref{coulomb}. When we move $m_2$ to position $P_2$ from infinity it will be acted upon the field created by $m_1$ and the work done will be
$$
\frac{m_1m_2}{r_{12}}.
$$
Now these two masses jointly create a new electrostatic field of joint magnitude
$$
\frac{m_1m_3}{r_{13}^2}+\frac{m_2m_3}{r_{23}^2}
$$
acting on any third mass $m_3$.  The work done  by moving $m_3$ to occupy
to position $P_3$ from infinity is then (by superposition)
$$
\frac{m_1m_3}{r_{13}} +\frac{m_2m_3}{r_{23}}.
$$
This now establishes yet a third electrostatic field whose magnitude of force will be given by
$$
\frac{m_1m_4}{r_{14}^2}+\frac{m_2m_4}{r_{24}^2}+\frac{m_3m_4}{r_{34}^2}
$$
and
that would act on a fourth particle $m_4$. Continuing the process for each of the $n$ masses by adding these all together will give the work done against the sequence of electrostatic fields
\begin{equation}\label{discrete}
E=\frac{1}{2}\sum_i\sum_j \frac{m_im_j}{r_{ij}},
\end{equation}
since we are identifying work and potential energy  (recall that work is a difference in potential energies by the law of conservation of energy) up to a fixed reference constant.  We write ``$E$'' here in place of ``$W$''.

Notice also that if we wrote equation \ref{discrete} in an iterated form
\begin{equation}
W=E=\frac{1}{2}\sum_im_i\sum_j \frac{m_j}{r_{ij}},
\end{equation}
the inside sum has an interesting physical interpretation. For a fixed position $P_i$
\begin{equation}\label{innerdiscrete}
U(i)=\sum_j \frac{m_j}{r_{ij}}
\end{equation}
represents the the potential energy (work) of bringing a unit charge (not necessarily $m_i$) to position $P_i$, when acted upon by the assemblage of masses $m_j$, $j\neq i$, which we denote by $U(i)$.

Our next step will be analogical. We ask what the continuous mass distribution (charge distribution) analogue of equation \ref{discrete} and \ref{innerdiscrete} will be?
Physically we imagine that a charge distribution consists of space filling charge as well as surface covering charge. For example, we might think of a conductor possessing charge in its interior as well as along its surface. Thus, we will consider both  three-dimensional and two-dimensional analogues.

Let us first consider the three-dimensional (volume) analogues of equations \ref{discrete} and \ref{innerdiscrete}. That is, we now move from a mass system in the plane to three-dimensional space.
 The analogy may be a little easier to see if we included the specific dependence of the masses $m_i$ on the points $P_i$ in equation \ref{discrete} \footnote{In these summations we assume $i\neq j$.} as
\begin{equation}\label{discreteya}
E=\frac{1}{2}\sum_i\sum_j \frac{m_i(P_i)m_j(P_j)}{|P_i-P_j|}=\frac{1}{2}\sum_im_i(P_i)\sum_j \frac{m_j(P_j)}{|P_i-P_j|}=\frac{1}{2}\sum_i m_i(P_i)U(P_i)
\end{equation} \label{star1}
where here the potential energy of equation \ref{innerdiscrete} is written
\begin{equation} \label{star2}
U(P_i)=\sum_j \frac{m_j(P_j)}{|P_i-P_j|}.
\end{equation}

If we go a little further and write the coordinates of $P_i$, $P_j$ as $(x_i,y_i,z_i)$ and
$P_j$ as $(r_j,s_j,t_j)$, the above can be expressed as

\begin{align*}
E=\frac{1}{2}\sum_i\sum_j \frac{m_i(x_i,y_i,z_i)m_j(r_j,s_j,t_j)}{\sqrt{(x_i-r_j)^2+(y_i-s_j)^2+(z_i-t_j)^2}}\\
=\frac{1}{2}\sum_im_i(x_i,y_i,z_j)\sum_j \frac{m_j(r_j,s_j,t_j)}{\sqrt{(x_i-r_j)^2+(y_i-s_j)^2+(z_i-t_j)^2}}.
\end{align*}
In this case, writing $U$ as a function of the coordinates of the $P_i$ amounts to saying
$$
U(x_i,y_i,z_i)=\sum_j \frac{m_j(r_j,s_j,t_j)}{\sqrt{(x_i-r_j)^2+(y_i-s_j)^2+(z_i-t_j)^2}}.
$$

Because these expressions are sums for a discrete mass distribution in  space, the continuous analogues will involve integrals over mass distributions in  space, in other words ``double-triple'' integrals. Replacing $m(x_i,y_i,z_i)$ with a charge density $\rho(x,y,z)$,
a natural continuous analogue  of  equation \ref{discreteya} and \ref{star2} is given by  the double-triple (six-fold integral)

\begin{align*}
E=\frac{1}{2}\iiint\iiint \frac{\rho(x,y,z)\rho(r,s,t)drdsdtdxdy}{\sqrt{(x-r)^2+(y-s)^2+(z-t)^2}}\\
=\frac{1}{2}\iiint \rho(x,y,z)\left[\iiint \frac{\rho(r,s,t)drdsdt}{\sqrt{(x-r)^2+(y-s)^2+(z-t)^2}}\right]dxdydz
\end{align*}
so that

\begin{eqnarray} \label{continuouspace}
E=\frac{1}{2}\iiint \rho(x,y,z) U(x,y,z)dydxdz, \quad\text{ where }\\  U(x,y,z)=\iiint \frac{\rho(r,s,t)drdsdt}{\sqrt{(x-r)^2+(y-s)^2+(z-t)^2}}.\nonumber
\end{eqnarray}
The above integrals are understood to be over all of space.\footnote{Note that $E$ is a spatial energy due to the spatial component of the charge distribution while $U$ is a corresponding spatial potential. A less sightly but possibly more accurate way of  writing these would be to denote them with subscripts like $E_{\mathbb{R}^3}$ an $U_{\mathbb{R}^3}$.} The  astute reader will notice that this last integral is the "Newtonian potential". Notice also that the volume potential $U(x,y,z)$ is defined at all points of space including on the surface $S$.

Next let us consider the two-dimensional surface analogues of equations \ref{discrete} and \ref{innerdiscrete}. We proceed along similar lines.

Returning briefly to the plane we now write the coordinates of $P_i$ as $(x_i,y_i)$ and
$P_j$ as $(r_j,s_j)$. The above can then be expressed as

\begin{eqnarray}\label{star3}
E=\frac{1}{2}\sum_i\sum_j \frac{m_i(x_i,y_i)m_j(r_j,s_j)}{\sqrt{(x_i-r_j)^2+(y_i-s_j)^2}} \nonumber \\
=\frac{1}{2}\sum_im_i(x_i,y_i)\sum_j \frac{m_j(r_j,s_j)}{\sqrt{(x_i-r_j)^2+(y_i-s_j)^2}}.
\end{eqnarray}
In this case, writing $U$ as a function of the coordinates of the $P_i$ amounts to saying
\begin{equation} \label{asterisk3}
U(x_i,y_i)=\sum_j \frac{m_j(r_j,s_j)}{\sqrt{(x_i-r_j)^2+(y_i-s_j)^2}}.
\end{equation}

Because these expressions are sums for a discrete mass distribution in the plane, the continuous analogues will be integrals over mass distributions in the plane, in other words ``double-double integrals''.  Replacing $m(x_i,y_i)$ with a charge density $\rho(x,y)$,
a natural continuous analogue of equation \ref{star3}  is given by  the double-double (four-fold integral)

$$
E=\frac{1}{2}\iint\iint \frac{\rho(x,y)\rho(r,s)drdsdxdy}{\sqrt{(x-r)^2+(y-s)^2}}
$$
$$
=\frac{1}{2}\iint \rho(x,y)\left[\iint \frac{\rho(r,s)drds}{\sqrt{(x-r)^2+(y-s)^2}}\right]dxdy=
\frac{1}{2}\iint \rho(x,y)U(x,y)dxdy,
$$
where
$$
U(x,y)=\iint \frac{\rho(r,s)drds}{\sqrt{(x-r)^2+(y-s)^2}}.
$$

The above reasoning also works if we are considering assemblages of   continuous electrical mass distributions on surfaces. Let us therefore make yet a further analogical leap and postulate that the surface analogue of equations  \ref{star3} and
\ref{asterisk3}  as being given by

\begin{eqnarray}\label{continuoussurface}
E=\frac{1}{2}\iint_S \rho(x,y)U(x,y)dS, \nonumber \\
 U(x,y)=\iint_S \frac{\rho(r,s)dS}{\sqrt{(x-r)^2+(y-s)^2+(f(x,y)-f(r,s))^2}}. 
\end{eqnarray}

These formulae need a little explanation. Here $dS$ represents (Lebesgue) surface area measure and coordinates $(x,y)$ are understood to be local coordinates on the surface. For example, if our surface in space were defined by an equation  $z=f(x,y)$, $dS$ could be replaced by a formula in multivariable calculus like $dS=\sqrt{1+|\nabla f(x,y)|^2}dxdy$. Since our interest is more conceptual than computational, we needn't make use of this. But we could if we wanted to write equation \ref{continuoussurface}
\begin{equation}\label{localcoords}
E=\frac{1}{2}\iint_{\mathcal{R}}  \rho(x,y)U(x,y)\sqrt{1+|\nabla f(x,y)|^2}dxdy
\end{equation}
where $\mathcal{R}$ is an appropriate parameter domain in the plane.

Likewise
\begin{equation}\label{Uxy}
U(x,y)=\iint_\mathcal{R} \frac{\rho(r,s)\sqrt{1+|\nabla f(r,s)|^2}}{\sqrt{(x-r)^2+(y-s)^2+(f(x,y)-f(r,s))^2}}drds
\end{equation}

Again, the points $(x,y)$  and $(r,s)$ are understood to be local coordinates on $S$ of points $\xi\in S$ and $\zeta\in S$. Thus, the potential energy function $U$ in equation \ref{continuoussurface} could be more succinctly expressed by
\begin{equation}\label{surfacepotential}
U(\xi)=\iint_S \frac{\rho(\zeta)}{|\xi-\zeta|}dS(\zeta),\qquad \xi\in S.
\end{equation}
Notice that although $U$ is initially defined on $S$, this integral is actually well defined all of
$\mathbb{R}^3$, and that it has the same physical interpretation of the potential energy at $\xi\in\mathbb{R}^3$ of a unit charge at $\xi\in\mathbb{R}^3$ relative to the charge distribution on $S$.  

With the understanding that $U$ is defined on all of $\mathbb{R}^3$,
 the equations for surface energy \ref{continuoussurface} and potential \ref{surfacepotential} yield
\footnote{As with the energy and potential determined by the volume charge distribution above, the energy and potential appearing in equations \ref{surfacepotential} and \ref{surfaceenergy} would be more accurately denoted by $E_S$ and $U_S$ to emphasize that this $E$ and $U$ are the surface energy  and the surface potential  of the surface charge distribution. The  ``total potential''  would then be defined as  $U=U_{\mathbb{R}^3}+U_S$. The total energy $E$, however, is not simply
$E_{\mathbb{R}^3}+E_S$, but requires an additional term for the mutual energy. This is a consequence of the functional defined by equation \ref{ei} being nonlinear. We address this shortly.}
\begin{equation} \label{surfaceenergy}
E=\frac{1}{2}\iint_S \rho(\xi)U(\xi)dS(\xi).
\end{equation}
This is what Kellogg \cite{Kel} called the {\it self-potential} of the charge distribution on $S$. When the charge distribution is in equilibrium we expect $E$ to be in some sense minimized. In \cite{Gr} Grube uses the German word {\it Wirkung} (action) to refer to the same concept, possibly reminiscent of Maupertuis' now obsolete terms for potential energy.

\section{Calculating the Dirichlet Energy} \label{dcalc}

\begin{quotation}
{\it The existence of a unique solution is very plausible by the "physical argument": any charge distribution on the boundary should, by the laws of electrostatics, determine an electrical potential as solution.}
-Hans Freudenthal
\end{quotation}

A fact known to Dirichlet was that if $\rho(\xi)$ is a continuous charge surface density on $S$, it can be recovered from the normal derivatives of the surface potential $U$ at the boundary:
\begin{equation}  \label{poissonsurface}
\rho(\xi)=-\frac{1}{4\pi}\left[ \frac{\partial U}{\partial \mathbf{n}_+}-\frac{\partial U}{\partial \mathbf{n}_-}\right].
\end{equation}
Here $\frac{\partial U}{\partial \mathbf{n}_+}$ is the radial (along the normal line) limit of the normal derivative from the outside of $S$ while $\frac{\partial U}{\partial \mathbf{n}_-}$ is the radial limit of the normal derivative from the inside of $S$.
This seems to have been originally discovered by Poisson \cite{Poi} and appears in George Green's famous essay \cite{Green}. See also \cite{Gr}, \cite{G}, and \cite{Kel}. The equation (\ref{poissonsurface}) can also be derived from Gauss' law.
Perhaps more familiar to many readers is Poisson's  volume density equation
\begin{equation}\label{poissonvolume}
\Delta U(x,y,z)+4\pi\rho(x,y,z)=0,
\end{equation}
where $U$ is the spatial potential of the spatial charge density $\rho$. Here both equations \ref{poissonsurface} and \ref{poissonvolume} are expressed in Gaussian units, though this is inessential to our purpose.

In summary, imagine that we have a closed bounded region $\Omega$ in space with boundary $S=\partial\Omega$. We imagine now space filled by a continuous charge distribution and that $S$ itself is covered by a lamina of charge distribution.  Separately, each part of the system, the volume distribution and the surface distribution has an energy, its self-potential. At first glance one might think  that the total energy or self-potential of the system as a whole would simply be the sum of the self-potentials of the parts. However, there will be a third source of energy that arises from the interaction of the two parts called the  mutual energy. To see why, let us provisionally assume that the self-potential (energy) of a charge distribution having total potential $U$ is written as an inner product.  In the case of a spatial or surface  charge distribution in isolation, these energies would take the forms 
$$
E_{\mathbb{R}^3}=\iiint |\nabla U_{\mathbb{R}^3}|^2dV \qquad \text{and}\qquad E_S=\iiint |\nabla U_S|^2dV. 
$$
At the same time we expect the total energy of the total potential $U=U_{\mathbb{R}^3}+U_S$ to likewise be given by equation \ref{ei}. Expanding the energy integral as an inner product then gives
$$
E=\iiint |\nabla \left( U_{\mathbb{R}^3}+U_S\right)|^2dV =\iiint |\nabla U_{\mathbb{R}^3}|^2dV +\iiint |\nabla U_S|^2dV+2\iiint \nabla U_{\mathbb{R}^3}\cdot \nabla U_S\,dV.
$$
One-half of this last integral,
$$
\iiint \nabla U_{\mathbb{R}^3}\cdot \nabla U_S\,dV,
$$
represents the mutual energy.\footnote{Called mutual potential by Kellogg.} It can be shown that the mutual energy can be expressed by either of the (necessarily equal) integrals
$$
\iiint \rho(x,y,z)U_S(x,y,z)\,dV\qquad\text{or}\qquad \iint_S \rho(x,y)U_{\mathbb{R}^3}(x,y)\, dS.
$$
Physically, the mutual energy represents the work done when the spatial and surface charge distributions are brought together from an infinite distance without changing form. See Kellogg \cite{Kel}.
The total energy (self-potential) is then obtained by adding together the equations for spatial and surface energies \ref{continuouspace}, \ref{surfaceenergy}  and either of the integrals for mutual energy. 
Kellogg \cite{Kel},  describes it thus:
\begin{quotation}
\dots if two bodies are brought, without change of form, from an infinite distance apart to a given position, the work done, or their mutual potential, is the integral over either body of the product of its density by the potential of the other. The self-potential of the system of the two bodies is the sum of the self-potential of the bodies separately and their mutual potential.
\end{quotation}

Henceforth, given densities $\rho(x,y,z)$ and $\rho(x,y)$, $U$ will denote the total potential, that is, the sum of the volume and surface potentials 
$$
U=U_{\mathbb{R}^3}+U_S.
$$ 
Let us further assume the spatial and surface densities $\rho(x,y,z)$ and $\rho(x,y)$ are ''sufficiently regular'' so that the potentials exist and sufficiently smooth. See \cite{Kel} Theorem III, chapter VI or volume II of Courant and Hilbert  \cite{CH} chapter IV page 246 for precise statements.
With this understanding, we claim the total energy can be written as 
 
\begin{equation} \label{energysummary}
E=\frac12\iiint_{\mathbb{R}^3} \rho(x,y,z)U(x,y,z)dxdydz + \frac12\iint_S \rho(x,y)U(x,y)dS.
\end{equation}

To see why, note that if we write $U=U_{\mathbb{R}^3}+U_S$, then equation \ref{energysummary} would take the form 
$$
\frac12\iiint_{\mathbb{R}^3} \rho(x,y,z)[U_{\mathbb{R}^3}+ U_S]dxdydz + \frac12\iint_S \rho(x,y)[U_{\mathbb{R}^3}+ U_S]dS
$$
$$
=\frac12\iiint_{\mathbb{R}^3} \rho(x,y,z)U_{\mathbb{R}^3}\, dxdyxz+ \frac12\iiint_{\mathbb{R}^3} \rho(x,y,z)U_S\,dxdydz + \frac12\iint_S \rho(x,y)U_{\mathbb{R}^3}\, dS+ \frac12\iint_S \rho(x,y)U_S\,dS.
$$
The first integral in equation \ref{energysummary} gives the self-potential of the spatial distributions when $U=U_{\mathbb{R}^3}$ while the second integral gives the self-potential of the surface distribution when $U=U_S$. When $U$ is taken to be the sum of the spatial and surface potentials, the same equation implicitly includes the mutual energy. These appear as the last two integrals.

\begin{remark}
Notice that  although the volume potential
\ref{continuouspace} and surface potential \ref{surfacepotential} both involve the same ``Newton kernel'' 
$$
\frac{1}{\sqrt{(x-\cdot)^2+(y-\cdot)^2+(z-\cdot)^2}},
$$
harmonic in the variables $x$, $y$, and $z$.  At points of space off $S$ the volume potential involves integration over singularities while the surface potential, which is over $S$, does not. Consequently, the surface potential will be harmonic, but by the Poisson formula \ref{poissonvolume} the volume potential will only be harmonic when it is zero corresponding to $\rho\equiv 0$. In essence the ``total potential'' =``volume potential'' + ``surface potential''= ``non-harmonic'' potential (except when 0) + ``harmonic'' potential. The total potential is harmonic only when the volume potential vanishes.
\end{remark}

The trick to relating the energy in equation \ref{energysummary} to the Dirichlet energy integral  \ref{ei} is to use Green's first identity\footnote{In \cite{Kel} Kellogg either misstates or misprints this as Green's second identity.}:
\begin{equation}\label{GFI}
\iint \left(A\,\frac{\partial B}{\partial\mathbf{n}}\right)\,dS=\iiint (A\,\Delta B)\, dV+
\iiint (\nabla A\cdot\nabla B)\, dV.
\end{equation}

\begin{remark}\label{remark1}
Here $\mathbf{n}$ is the usual outward unit normal vector to the surface. 
The expression $$\frac{\partial B}{\partial\mathbf{n}}$$ is the outward pointing normal derivative in $\Omega$  identified with 
$\frac{\partial B}{\partial\mathbf{n}_-}$ while $\frac{\partial B}{\partial\mathbf{n}_+}$ is identified with the outward pointing normal derivative for the region
 $ext(\Omega)$. Intuitively we may think of this as the normal to $S$ pointing into the interior of $\Omega$, so that $-\frac{\partial B}{\partial\mathbf{n}_+}$ also points outward from $\Omega$.
\end{remark}

Given a charge distribution in space, it will have a potential function $U$. Note that $U$ need not be harmonic. Conversely, suppose we begin with a sufficiently smooth function $U(x,y,z)$ defined on all of space. The Poisson density formulas
\ref{poissonsurface} and \ref{poissonvolume} then determine spatial and surface charge distributions.
More precisely, let
 $U:\mathbb{R}^3\to\mathbb{R}^1$. We assume that $U$ is continuous up to $S$ and twice continuously differentiable on $\Omega\cup ext(\Omega)$. Motivated by the Poisson density formulas
  \ref{poissonvolume} and \ref{poissonsurface}, we define a volume density
\begin{equation} \label{pd}
\rho_{\mathbb{R}^3}:=-\frac{1}{4\pi}\Delta U,
\end{equation}
and a surface  density
\begin{equation} \label{sd}
\rho_S=-\frac{1}{4\pi}\left[ \frac{\partial U}{\partial \mathbf{n}_+}-\frac{\partial U}{\partial \mathbf{n}_-}\right].
\end{equation}
Although it requires a mathematical proof, it is physically plausible that 
 \begin{equation}\label{total}
U=U_{\mathbb{R}^3}+U_S
\end{equation}
where $U_{\mathbb{R}^3}$ and $U_S$ are defined to be the volume and surface potentials $\rho_{\mathbb{R}^3}$ and $\rho_S$ given in equations \ref{pd} and \ref{sd}. We also assume that $U$ is sufficiently regular at infinity so that these potentials both exist.

Substituting these densities along with the given function $U$ expressed as a total potential as in equation \ref{total} into equation \ref{energysummary}
then yields
\begin{equation} \label{E}
E=\frac{-1}{8\pi}\left[ \iiint_{\mathbf{R}^3} U\Delta U\, dV+\iint_S U\left(\frac{\partial U}{\partial \mathbf{n}_+}-\frac{\partial U}{\partial \mathbf{n}_-}\right)dS\right]
\end{equation}
Unpacking the volume and surface integral in \ref{E} we see
\begin{eqnarray}
\iiint_{\mathbf{R}^3} U\Delta U\, dV=\iiint_\Omega U\Delta U\, dV + \iiint_{ext(\Omega)} U\Delta U\, dV,\label{32}
\\
\iint_S U\left(\frac{\partial U}{\partial \mathbf{n}_+}-\frac{\partial U}{\partial \mathbf{n}_-}\right)dS=\iint_S U\frac{\partial U}{\partial \mathbf{n}_+}dS-\iint_S U\frac{\partial U}{\partial \mathbf{n}_-}dS.  \label{33}
\end{eqnarray}

Applying Green's first identity  \ref{GFI} to each integral on the right side of equation
\ref{33} in light of remark \ref{remark1} yields
\begin{eqnarray}
\iint_S U\frac{\partial U}{\partial \mathbf{n}_+}dS=-\iint_S U\left[-\frac{\partial U}{\partial \mathbf{n}_+}\right]dS \nonumber \\
=-\left(\iiint_{ext(\Omega)}U\Delta U\, dV+\iiint_{ext(\Omega)} |\nabla U|^2dV \right),\label{license}
\end{eqnarray}
and
$$
\iint_S U\frac{\partial U}{\partial \mathbf{n}_-}dS=\iiint_{\Omega} U\Delta U\, dV + \iiint_\Omega |\nabla U|^2 dV.
$$
Notice that the application of Green's first identity in equation \ref{license} is analogical. Strictly speaking, Green's identities apply to closed bounded regions. In the spirit of analogical license, we reason as if it also applied to the unbounded region
 $ext(\Omega)$.

Thus we rewrite equation \ref{33} as
\begin{multline*}
\iint_S U\left(\frac{\partial U}{\partial \mathbf{n}_+}-\frac{\partial U}{\partial \mathbf{n}_-}\right)dS=\\
-\iiint_{ext(\Omega)}U\Delta U\, dV-\iiint_{ext(\Omega)} |\nabla U|^2dV  -\iiint_{\Omega} U\Delta U\, dV - \iiint_\Omega |\nabla U|^2 dV.
\end{multline*}
Adding this to equation \ref{32} finally gives
\begin{eqnarray} 
E=\frac{-1}{8\pi}\left[- \iiint_ \Omega|\nabla U|^2\, dV-\iiint_ {ext(\Omega)}|\nabla U|^2\, dV\right] \nonumber \\
=\frac{1}{8\pi}\iiint_ \Omega|\nabla U|^2\, dV
+\frac{1}{8\pi}\iiint_ {ext(\Omega)}|\nabla U|^2\, dV. \label{complete-energy}
\end{eqnarray}

The modern incarnation of the notions of volume and surface potential arising from volume and surface densities described above corresponds to the notion of the potential of a finite Borel measure. What Kellogg in 1929 dubbed self-potential and mutual potential corresponds to the Coulomb energy of a measure and the mutual energy of a pair of measures. A good exposition of these ideas is contained in chapter three of Barry Simon's book \cite{simon3}.

\section{Electrostatic Interpretation of Minimum Dirichlet Energy}\label{eim}

\begin{quotation} 
{\it Similary, in 1847, Sir William Thomson, Lord Kelvin, attempted to found a proof on the least value of an integral. The same considerations were used by Dirichlet in lectures during the following decade $\cdots$

One might be led to it as follows. We imagine the region $R$, for which the problem is to be solved, and the rest of space, filled with charges, and in addition, a spread on the bounding surface $S$. -Oliver Dimon Kellogg
}
\end{quotation}

It is a basic fact of physics that in an electrostatic field, there can be no free charges (in stable equilibrium) in the interior of a conductor; all such charge accumulates on the surface of the conductor. One may imagine a negative charge distribution in the interior of 
$\Omega$ in which the repulsive force between like charges impells the charges to attempt to ``escape'' by congregating on the bounding surface. For point charges this is sometimes called Earnshaw's theorem after Samuel Earnshaw (1805-1888). See also the Feynman lectures 5.2 for a more detailed explanation in \cite{Feyn}.

Now let us suppose $f:S\to\mathbb{R}$ be a continuous function. It may be thought of as a prescribed boundary potential. As in definition \ref{admiss}, let $\mathcal{A}_f$ be the class of admissible extensions of $f$. Each extension $U$ may be regarded as a potential to which there is a corresponding charge distribution determined by the Poisson formulae \ref{poissonsurface} and \ref{poissonvolume}.\footnote{Recall that $U$ is understood to be the sum of a volume and surface potential, i.e. $U=U_{\mathbb{R}^3} +U_S$.} The only way $U$ can be spatially harmonic is if its volume potential is zero. This can only occur if all charges have congregated on $S$. This will be the case, for example, if the charge distribution is in electrostatic (stable) equilibrium. By the principle of minimum potential energy, stable equilibrium occurs when the total energy $E$ in equation \ref{complete-energy} is minimum amongst all potentials $U$ subject to the constraint that $U\in\mathcal{A}_f$. These observations constitute the heuristics behind the electrostatic argument for the Dirichlet principle. We expect that $D(U)$ in equation \ref{ei} above should have a minimum value, and, because it would correspond to a stable equilibrium, $U$ would be the sought after harmonic extension to $\Omega$.
\footnote{Although above we defined the energy by integrals over the whole space $\mathbb{R}^3$, its minimization entails that the integrals over both $\Omega$ and $ext(\Omega)$ must  be separately minimized. We remark more on this below.}

Since this argument is the crux of this paper, it seems appropriate to give the argument in more precise mathematical notation.  Denote by $\mathcal{J}(U,\rho,F)$\footnote{We choose the symbol $\mathcal{J}$ here to denote the units of energy Joules.} the energy (self-potential), now given by equation \ref{complete-energy}, of the datum (electrical system) given in definition \ref{datum}. For a prescribed boundary potential $f:\partial \Omega\to\mathbb{R}$, define\footnote{Notice we are here identifying the class of admissible extensions $\mathcal{A}_f$ with a class of electrical systems.}
$$
\mathcal{A}_f=\left\{(U,\rho,F): U\vert_{\partial\Omega}=f \right\},
$$
and
$$
B=\left\{ \mathcal{J}(U,\rho,F):(U,\rho,F) \in \mathcal{A}_f \right\}.
$$
As a set of nonnegative real numbers, $B$ is bounded below by zero and thus has a non-negative greatest lower bound. Supposing it to be a minimum of $B$, let $(U,\rho,F)$ be a datum in $\mathcal{A}_f$ for which the potential energy
$\mathcal{J}(U,\rho,F)$ is minimum. By the principle of potential energy minimization, this datum $(U,\rho,F)$ represents an electrical system in stable equilibrium. Thus, its charge density must vanish in $\Omega\cup ext(\Omega)$ and be fully concentrated on the boundary $\partial \Omega$.  This implies:
\begin{enumerate}
\item $U$ is harmonic in $\Omega\cup ext(\Omega)$.
\item $\iiint_{\mathbf{R}^3} |\nabla U|^2 dV$ is minimized over $\mathcal{A}_f$.
\end{enumerate}

To summarize: Given prescribed boundary value $f$, the extension $U$ of minimum Dirichlet energy $D(U)$ will be harmonic and thus solve the Laplacian boundary value problem for boundary data $f$.

The foregoing can now be turned around to furnish a strategy for solving the Dirichlet principle as follows.  Let $\Omega$ be a bounded open region in $\mathbb{R}^3$ with surface $S=\partial\Omega$. Let $f:S\to\mathbb{R}^3$ be a given sufficiently smooth prescribed boundary function.  For an appropriate class of admissible functions
$U:\overline{\Omega}\to \mathbb{R}^3$, $\left. U\right\vert_{S}=f$, find such
$U$ which minimizes the ``energy''
\begin{equation} \label{strategy}
D(U)=\iiint_{\Omega} |\nabla U|^2 dV.
\end{equation}
Then show that the minimizing function $U$ is harmonic in $\Omega$. Note that we have replaced the region of integration $\mathbb{R}^3$ in integral  by $\Omega$. This is justified by the fact that the minimum of the integral over $\mathbb{R}^3$ entails the minimum of the integral \ref{strategy}.  (The same will also be true of the integral over $ext(\Omega)$, but we need not concern ourselves with this because our interest is in the interior Dirichlet problem.)

\section{The Dirichlet Principle as a Mathematical Problem in the Calculus of Variations}\label{cv}
Having arrived at a possible theorem by physical reasoning, we would like to solve the Dirichlet problem using the method of energy minimization suggested by physical considerations above  but do so in a mathematically rigorous manner.
  We present this as a ``contingent theorem'' instead of a theorem to highlight the fact that even though we are no longer reasoning  by physical analogy, the theorem still remains unproved without the further refinements of the hypothesis supplied later by Hilbert in \cite{H}.

\begin{cthm*} Suppose $\Omega\subset\mathbb{R}^3$ is a bounded open connected region bounded by the surface $S$.  Let $f:S\to\mathbb{R}$ be a given continuous function on $S$.  Then there exists a (unique) harmonic function $u:\overline{\Omega}\to\mathbb{R}$, for which $u=f$ on $S$.
\end{cthm*}

Our strategy, suggested by the above physical reasoning, will be to argue that for a given boundary value $f$, the extension $u$ that minimizes the nonlinear functional 
\ref{strategy} is the sought after harmonic extension.
Following Kellogg \cite{Kel}, we present the classical variational argument. 
To this end let us define the integral
\begin{equation} \label{duv}
D(u,v)=\iiint_\Omega \nabla u\cdot\nabla v\, dV
\end{equation}
for any two continuously differentiable functions $u$ and $v$. When $u=v$, we just write $D(u)$ as in equation \ref{strategy}.

If $B\subset \Omega$ is a sphere (open ball) having boundary $\partial B$, we will use the notation
$$
D_B(u)=\iiint_{B} |\nabla u|^2 dV,\quad D_B(u,v)=\iiint_B \nabla u\cdot\nabla v\, dV
$$
to denote the integrals corresponding to equations \ref{strategy} and \ref{duv} but over the set $B$.

Let us further define the class $\mathcal{\mathcal{A}}$ of functions $u$ twice continuously differentiable in $\Omega$ and continuous up to the boundary, that is on $\overline{\Omega}$. Given $f:S\to\mathbb{R}$, let $\mathcal{A}_f$  denote those $u\in\mathcal{A}$ such that $u=f$ on $S$. 
Clearly the conjectured theorem is true for constant $f$. So assume $f$ is not constant. For each admissible function $u\in \mathcal{A}_f$, it is clear that $D(u)\geq 0$. Thus the set of all such values of $D(u)$ is bounded below and so possess a positive greatest lower bound. (Remember we are assuming $f$ is not constant.) {\bf Let us assume this greatest lower bound is also a minimum.}  This is a priori not obvious. Infimums and minimums are not always the same. In fact it is not in general true here without more careful hypotheses. But let us assume that $u\in\mathcal{A}_f$ minimizes $D(u)$.

Observe that the minimizer $u$ must also minimize $D_B(u)$ over the extensions of
$u|_{\partial B}$. This is because since $u\in \mathcal{A}_f$, it's values could be altered inside $B$ to get a smaller minimum.

Next choose a function $\tilde{u}$  on the closed sphere $\overline{B}$ that is an admissible extension of $u$ on $\partial{B}$ to all of $B$. That is, $\tilde{u}|_{\partial{B}}=u|_{\partial B}$.

Define the function $h=u-\tilde{u}$. Note that on $\partial{B}$, $h=0$. For $x$ real we will make some observations about the function
$$
u+xh,
$$
which has boundary value  on $\partial B$ as $u$ itself. Hence $u+xh$ is an admissible extension of $u$ on the  boundary of the sphere to all of $B$.
We emphasize that $u$, $h$, and $u+xh$ are functions defined on the closed sphere, but $x$ is just a real number. 

 It's an easy calculation to show
\begin{equation} \label{Iexpand}
D_B(u+xh)=D_B(u)+2xD_B(u,h)+x^2D_B(u).
\end{equation}
By assumption $D_B(u)$ is minimal amongst functions in $\mathcal{A}_f$ so
$D_B(u)\leq D_B(u+xh)$, for all $x$. Thus $D_B(u+xh)-D_B(u)\geq 0$. Equation \ref{Iexpand} then implies that
\begin{equation} \label{nonzero}
2xD_B(u,h)+x^2D_B(u)\geq 0.
\end{equation}
Since inequality \ref{nonzero} is understood to hold for all real $x$, it is not possible that $D_B(u)=0$ while $D_B(u,h)\neq 0$. Suppose that $D_B(u)\neq 0$. 
Consider now the quadratic in $x$ defined by  $Q(x)=2xD_B(u,h)+x^2D_B(u)=xD_B(u) \left[x+\frac{2D_B(u,h)}{D_B(u)}\right]$,  By equation
\ref{nonzero}, $Q(x)\geq 0$ for all $x$, hence $D_B(u,h)=0$. By Green's first identity
\begin{equation}\label{afortiori}
0=D_B(u,h)=\iiint_{B}\nabla u\cdot\nabla h\, dV=\iint_{\partial B} h\frac{\partial u}{\partial\mathbf{n}}dS-\iiint_B h\Delta u\, dV.
\end{equation}
 Moreover, since $h=0$ on $\partial{B}$, it follows that
\begin{equation} \label{almost}
\iiint_B h\Delta u\, dV=0.
\end{equation}
Equation \ref{almost} is understood to hold for every function $h=u-\tilde{u}$ as described above. We claim this allows us to conclude that $\Delta u\equiv 0$ in $\Omega$. To see why suppose there were a point $p\in\Omega$ for which
$\Delta u(p)>0$. Let $B$ be a small sphere centered at $p$   contained in $\Omega$
on which $\Delta u> 0$.
  Choose  $\tilde{u}$  so that  $h>0$ inside the sphere.  (Note this is possible because there is no requirement that $\tilde{u}$ be harmonic.) But this would imply
$$
\iiint_B h\Delta u\, dV>0,
$$
contradicting equation \ref{almost}.  The possibility that $\Delta u(p)<0$ is disposed of in the same way. If $D_B(u,h)=0$ but $D_B(u)\neq 0$, the equations \ref{afortiori} and \ref{almost} hold a fortiori and the conclusion is the same.

Thus, we conclude that {\bf if} the Dirichlet energy integral has a minimizer, this minimizer is a harmonic function which solves the Dirichlet problem. The uniqueness of this solution follows, for example, from the maximum principle for harmonic functions.

\section{Epilogue}

In order to see the forest from the trees, let us give a brief recap of our chain of logical, physical, and analogical reasoning.

\begin{enumerate}
\item Given a prescribed boundary function $f$, we consider the set of all admissible extensions to $\Omega\  \cup\  ext(\Omega)$, and denote this by $\mathcal{A}_f$.
\item Each admissible extension $U$  will serve as a total($U=U_{\mathbb{R}^3}+U_S$) potential for a corresponding charge distribution with spatial and surface components.  These densities can be recovered from the Poisson formulae \ref{poissonsurface} and \ref{poissonvolume}.
\item Any potential $U$ will satisfy Laplace's equation (be harmonic) in $\Omega\ \cup\  ext(\Omega)$ precisely when its volume component vanishes.  (Its harmonicity comes from the surface potential.)
\item The volume potential will vanish when the spatial charge density is zero. (When there is no charge in the interior.) 
\item Physically we expect there to be a potential $U$ of minimum
energy
\begin{equation} \label{DE}
D(U)=\iiint_{\Omega} |\nabla U|^2 dV.
\end{equation}
Here we are ignoring the corresponding integral over $ext(\Omega)$ as both would have to be minimized.
\item This minimum energy potential should correspond to a stable equilibrium in which there can be no charge present in 
$\Omega\ \cup \ ext(\Omega)$. 
\item No charge present off the surface of the conductor implies the volume potential is zero, hence the total potential equals the surface potential which is harmonic. 
\item Because $U$ was defined to be an extension of the prescribed boundary value $f$, it constitutes a harmonic extension of $f$.
\end{enumerate}
The upshot of all of this is that to solve a Dirichlet problem with boundary data $f$, we should look for the extension $U$ that minimizes  the integral \ref{DE} and try to show that it is harmonic.

\section{Historical Notes}
Vague allusions to the folklore of  the Dirichlet principle are often found in popular accounts of the history of mathematics. See for example \cite{Reid} chapter IX ``Problems''.
Our presentation is based on, inspired by,  and indebted to that of Kellogg \cite{K} which is itself is based on Grube \cite{Gr}. The foregoing sections \ref{ho}, \ref{eim}, and \ref{cv}   are essentially reproductions of Kellogg with additional mathematical, physical, and historical details.
Monna's book \cite{M} contains additional
helpful commentary as well as an extensive bibliography on nineteenth century potential theory.  Less extensive but likewise additional helpful commentary and references are found in Kline \cite{MK}. See also Gorkin and Smith \cite{GS} for a broad overview of the quest to put the Dirichlet principle on firm footing that culminated with  Hilbert \cite{H}. Garding's article \cite{Gard} contains an overview of the Dirichlet problem from it's prehistory involving Green and Gauss to the period after Riemann
 and the impetus it provided to the later research of  Weierstrass, Neumann, Schwarz, Poincar\'{e}, Hilbert, Perron, and F. Riesz.

The minimization problem for the integral \ref{strategy}  appeared (without the energy interpretation) in George Green's 1835 paper 
\cite{Green1835}. A similar integral minimization problem was considerd in 1847 by William Thomson \cite{T}. Gauss' seminal 1839 paper \cite{Gau} was one of the foundational papers in potential theory. It contains physical reasoning based on electromagnetic theory and the principle of minimum potential energy for charge distributions similar to the foregoing. Geppert's commentary
\cite{Gep} is a helpful supplement as is Schaefer's \cite{S}. Bacharach's book \cite{Bach} contains an extensive  history of potential theory circa the late nineteenth century.

\end{document}